\documentclass[12pt]{article}
\usepackage{amsmath}
\usepackage{amsfonts}
\usepackage{latexsym}
\usepackage{amssymb}

\newcommand{\Z}{{\mathbf Z}}

\newcommand{\R}{{\mathbf R}}
\newcommand{\C}{{\mathbf C}}

\newcommand{\kk}{{\mathbf k}}
\newcommand{\cat}{{\rm {cat }}}
\newcommand{\id}{{\rm {id }}}

\newcommand{\Hom}{{\rm Hom}}
\newcommand{\im}{{\rm im}}

\newcommand{\grad}{{\rm grad}}

\newtheorem{theorem}{Theorem}[section]

\newtheorem{lemma}[theorem]{Lemma}

\newtheorem{example}[theorem]{Example}
\newtheorem{definition}[theorem]{Definition}
\numberwithin{equation}{section}

\newenvironment{proof}{\par{\bf Proof.} }{\par}

\title{Zeros of closed 1-forms, homoclinic orbits, \\
and Lusternik - Schnirelman theory}
\author{Michael Farber\footnote{Partially supported by the US - 
Israel Binational
Science Foundation; part of
this work was
done while the author visited Max-Planck Institute for Mathematics in 
Bonn}}
\date{June 18, 2001}

\begin{document}

\sloppy

\maketitle

\begin{abstract}
In this paper we study topological lower bounds on the number of zeros of closed 1-forms
without Morse type assumptions. We prove that one may 
always find a representing closed 1-form having at most one zero. We 
introduce and study a generalization $\cat(X,\xi)$ of the notion of Lusternik - Schnirelman category, 
depending on a topological space $X$ and a 1-dimensional real cohomology class $\xi\in H^1(X;\R)$.
We prove that any closed 1-form $\omega$ in class $\xi$
has at least $\cat(X,\xi)$ zeros assuming that $\omega$ admits a gradient-like vector field with no homoclinic cycles.
We show that the number $\cat(X,\xi)$ can be estimated from below in terms of the cup-products and higher
Massey products.

This paper corrects some my statements made in \cite{F0}, \cite{F1}.


\end{abstract}

\section{Introduction}
The Novikov inequalities \cite{N, N1} estimate the numbers of zeros of
different indices of closed 1-forms $\omega$ on manifolds lying in a given cohomology
 class, assuming that all the
zeros are non-degenerate in the sense of Morse. 
Applications of the Novikov inequalities in mechanics, in geometry and in symplectic topology are well-known, see 
\cite{EG},  \cite{HS}, \cite{S1}, \cite{S2}, \cite{VO}. 

In this paper we show that one may always realize a nonzero cohomology 
class by a closed 1-form with at most one (degenerate) zero.
In the proof we use the technique of rearrangements of critical points and the result of F. Takens \cite{T},
which describes the conditions when several critical points of a function could be collided into one. 

Central role in this paper plays a suitable generalization of the notion of Lusternik - Schnirelman category. 
For any pair $(X,\xi)$, consisting of a topological space $X$ and a real cohomology class $\xi\in H^1(X;\R)$,
we define a non-negative integer $\cat(X,\xi)$, the category of $X$ with respect to the cohomology class
$\xi$. The definition of
$\cat(X,\xi)$ is similar in the spirit to the definition of $\cat(X)$; it deals with open covers of $X$ with certain homotopy
properties. 
We show that $\cat(X,\xi)$ depends only on the homotopy type of $(X,\xi)$ and coincides with 
$\cat(X)$ in the case $\xi=0$. If $\xi\neq 0$ then $\cat(X,\xi) < \cat(X)$; we show by examples that the difference
$\cat(X) - \cat(X,\xi)$ may be arbitrary. 

The main Theorem of the paper  (Theorem \ref{lsmain1}) 
states that any smooth closed 1-form $\omega$ on a smooth closed 
manifold $X$ must have at least $\cat(X,\xi)$ geometrically distinct zeros, where
$\xi=[\omega]\in H^1(X;\R)$ denotes the cohomology class of $\omega$,
assuming that $\omega$ admits a gradient-like vector field with 
no homoclinic cycles. Recall that {\it a homoclinic orbit} is defined as a trajectory $\gamma(t)$, $t\in \R$,
 such that both limits 
$\lim_{t\to + \infty}\gamma(t)$ and $\lim_{t\to - \infty}\gamma(t)$
exist and are equal. More generally, {\it a homoclinic cycle of length $n$}
is a sequence of orbits $\gamma_1, \gamma_2, \dots, \gamma_n$
such that 
$$\lim_{t\to +\infty} \gamma_i(t) = \lim_{t\to -\infty} \gamma_{i+1}(t)$$
for $i=1, \dots, n-1$ and 
$$\lim_{t\to +\infty} \gamma_n(t) = \lim_{t\to -\infty} \gamma_{1}(t).$$ 

Viewed differently, the main Theorem of the paper claims that
any gradient-like vector field of a closed 1-form has a homoclinic cycle assuming that the number of zeros 
is less than $\cat(M,\xi)$. 

The homoclinic orbits were discovered by H. Poincar\'e and were studied by S. Smale. 
In the mathematical literature there are many results 
about existence of homoclinic orbits in Hamiltonian systems.
Homoclinic orbits may not exist in the gradient systems for functions, 
i.e. in the case $\xi=0$, corresponding to the 
classical Lusternik - Schnirelman theory. 

In recent papers \cite{F}, \cite{F0}, \cite{F1} we described cohomological 
cup-length type estimates on the number of zeros of closed 1-forms. Unfortunately, they are incorrect as stated.
 Purely algebraic Proposition 3 of \cite{F0} and a similar in character Lemma 6.6 of \cite{F1} are incorrect.
These algebraic statements hold under slightly stronger assumptions, however, in oder to meet these assumptions
in the main theorems, one has to make extra assumptions on closed 1-form. 
This paper gives a different (more geometric) approach to show that the main results of 
\cite{F1} hold for closed 1-forms having gradient-like vector field with no homoclinic cycles. 
The cohomological lower bounds 
of \cite{F0} require also some additional
assumption, which will be discussed elsewhere.

I would like to thank Pierre Milman, Kaoru Ono and Shmuel Weinberger for useful discussions and help.

\section{Colliding the critical points}

In this section we will prove the following realization result.

\begin{theorem} \label{lsreal}
Let $M$ be a closed connected $n$-dimensional smooth manifold, and let 
$\xi\in H^1(M;\Z)$ be a nonzero cohomology class. Then there exists a smooth closed 1-form $\omega$ in 
class $\xi$ having at most one zero.
\end{theorem}

\begin{proof} We will assume that the class $\xi\in H^1(M;\Z)$ is indivisible, i.e. is not a
multiple of another integral class. Our statement clearly follows from this special case.

Our purpose is to show that we may find a smooth map $\phi: M\to S^1$
with Morse critical points, having the following properties:

 (A) {\it The cohomology class of the closed 1-form 
$\tilde \omega \, = \, \phi^\ast(d\theta)$ coincides with $\xi$, where
$d\theta$ denotes the standard angular form on the circle $S^1$. }

(B) {\it All fibers $\phi^{-1}(b)$ are connected, where $b\in S^1$.}

(C) {\it The map $\phi$ has at most one critical value $b_0\in S^1$. In other words, all critical points of 
$\phi$ lie in the same fiber $f^{-1}(b_0)$. }

Having achieved this, we may apply the technique of F. Takens \cite{T}, pages 203 - 206,
which allows to collide the critical points 
of $\phi$ (equivalently, the zeros of closed 1-form $\tilde \omega$)
into a single degenerate critical point of a closed 1-form $\omega$ lying in the same cohomology class. 
Namely, firstly, we may find a piecewise smooth tree $\Gamma\subset \phi^{-1}(b_0)$ containing all the critical points
of $\phi$. Secondly, we may find a continuous map $\Psi: M\to M$ with the following properties:

{\it $\Psi(\Gamma)$ is a single point $p\in \Gamma$;

$\Psi|_{M-\Gamma}$ is a diffeomorphism onto $M-p$;

$\Psi$ is the identity map on the complement of a small neighborhood of $\Gamma$;
in particular, $\Psi$ is homotopic to the identity map $M\to M$.}

\noindent
The circle valued map $\phi\circ \Psi^{-1}$ is well-defined and is continuous. Moreover, 
$\phi\circ \Psi^{-1}$ is smooth on $M-p$.
Applying Theorem 2.7 from \cite{T} we see that we may replace the map $\phi\circ \Psi^{-1}$ in 
a small neighborhood of $p$ by a smooth map $\psi: M\to S^1$ having a single critical point. 
$\psi$ is homotopic to $\phi$ and thus closed 1-form 
$$\omega\, =\, \psi^\ast(d\theta)$$
 lies in the cohomology class
$\xi$ and has possibly a single zero.

In the rest of the proof 
we will show that we may find a smooth Morse map $\psi:M\to S^1$ with properties (A) - (C) above. 

It is well known that any indivisible class $\xi\neq 0$ 
may be realized by a connected codimension one submanifold $V\subset M$ 
with oriented normal bundle. Cutting $M$ along $V$ produces a compact
cobordism $N$ with $\partial N=\partial_+N\cup \partial_-N$, a disjoint
union of two copies of $V$. Consider a Morse function $f: N\to [0,1]$ 
having $0$ and $1$ as regular values and $f^{-1}(0)=\partial_+N$, $f^{-1}(1)=\partial_-N$.
We may assume that $f$ has no critical points of indices $0$ and $n=\dim M$. 
Moreover, we may construct $f$, such that all level sets $f^{-1}(c)$, where $c\in [0,1]$,
are connected and having 
the self-indexing property: all critical points of $f$ having Morse index $i$ lie in $f^{-1}(i/n)$, where 
$i=1, \dots, n-1$. The map $N\to S^1$, where $x\mapsto \exp(2\pi i f(x))$, defines a smooth map
$\phi_1: M\to S^1$ in the cohomology class $\xi$ with connected fibers having the following property: 
for any critical point $m\in M$, $d{\phi_1}(m)=0$, 
with Morse index $i$, the image $\phi_1(x)$ equals $\exp(2\pi i/n)$. In other words, the critical points
of $\phi_1$ with the same Morse index lie in the same fiber, and these critical fibers $\phi_1^{-1}(b)$ appear
in the order of their Morse indices, while the point 
$b$ moves in the positive direction along the circle $S^1$.

For points $b_1, \dots, b_{r}\in S^1$ on the circle $S^1$ we will write 
$$b_1<b_2<\dots <b_{r}<b_1$$
to denote that moving from point $b_1$ in the positive direction along the circle $S^1$, we first meet
$b_2$, then $b_3, \dots$, until we again meet  $b_1$.

Let us now formulate the following approximation to property (C):

 (C$_j)$ {\it All critical points of a smooth Morse map 
$\phi: M\to S^1$ with Morse index $i$ lie in the same fiber $\phi^{-1}(b_i)$, where
$b_i\in S^1$, $i=1, \dots, n-1$, and
$$b_1<b_1<\dots <b_j=b_{j+1}=\dots=b_{n-1}<b_1.$$
In particular, 
the critical values $b_j=b_{j+1}=\dots = b_{n-1}$ coincide.
}

Note that (C$_1$) is equivalent to (C), which is our purpose.

We have found above a smooth map $\phi_1:M\to S^1$, satisfying (A), (B), and (C$_{n-1}$). 
On the next step we will show that we may replace $\phi_1$ by a smooth Morse map $\phi_2: M\to S^1$
with properties (A), (B) and (C$_{n-2}$). Let $b_1<b_2<\dots<b_{n-1}<b_1$ be the critical values of $\phi_1$.
Consider a point $c\in S^1$, lying between $b_{n-2}$ and $b_{n-1}$. Cut $M$ along the submanifold $\phi_1^{-1}(c)$
and consider the obtained cobordism $N$ and a Morse function $g$ from $N$ to an interval, 
obtained by cutting the circle $S^1$
at point $c$.
All level sets of $g$ are connected.
Moving from the
bottom of this cobordism to the top,
we meet $n-1$ critical levels; first we meet the level containing the critical points of Morse index $n-1$, 
then the levels containing the critical points with Morse indices $1, 2,\dots, n-2$. 
We will use the theory of S. Smale
of rearrangement of critical points. 
Choose a generic gradient like vector field $v$ for $g$. Then 
all integral trajectories of $\pm v$, which go out of the critical points of index $n-1$, reach $\partial N$ without
interaction with the other critical points.
Therefore we may slide the critical points of index $n-1$ some distance up, putting them on the same level
with the critical point of index $n-2$, see \cite{M}, Theorem 4.1. In other words, we may
replace $g$ by a Morse function $g'$, which 
coincides with $g$ near $\partial N$, has the same critical points, but the value at the critical points of index $n-1$ 
equals the value at the critical points of index $n-2$. Note that the level sets
of $g'$ are all connected:  

(a) the bottom level ${g'}^{-1}(0)$ is unchanged and so it is connected;

(b) passing the critical levels with Morse indices
$1, 2, \dots, n-3$ may not create nonconnected level sets;

(c) in principle, nonconnected level sets may appear after passing the top critical value contained the critical points of indices
$n-2$ and $n-1$; however in our situation
all higher upper level sets are the same as for the previous function $g$, and so they are all connected.

Folding this cobordism back, gives a smooth map $\phi_2: M\to S^1$
having properties (A), (B), and (C$_{n-2}$).

We may proceed similarly to find a smooth map $\phi_3: M\to S^1$ with properties (A), (B), and (C$_{n-3}$).
The critical values of $\phi_2$ are $b_1<b_2< \dots <b_{n-2}=b_{n-1}<b_1$.
We find a point $c\in S^1$ between $b_{n-3}$ and $b_{n-2}=b_{n-1}$ and cut $M$ along 
$\phi_2^{-1}(c)$. The Morse function on the obtained cobordism will have $n-2$ critical levels.
The lowest will be the level containing all critical points of indices $n-2$ and $n-1$ and then the critical levels 
of points with the Morse indices $1, \dots, n-3$. Repeating the above procedure, we may slide the critical points 
of index $n-2$ and $n-3$ the same distance up putting them on the same level with the critical points of index $n-3$.

Proceeding in this way inductively we arrive at a smooth map $\phi_{n-1}:M\to S^1$ 
having properties (A), (B), and (C$_{1}$) = (C).

This completes the proof. $\Box$
\end{proof}

According to a remark in \cite{EG}, a statement in the spirit of Theorem \ref{lsreal}
was made by Yu. Chekanov at a seminar talk 
in 1996. As far as I know, no written account
of his work is available.

\section{Category of a space with respect to a cohomology class}
\subsection{Definition of ${\rm cat}(X,\xi)$}

Let $X$ be a finite CW-complex and let $\xi\in H^1(X;\R)$ be a real cohomology class. 
We will define below a numerical invariant ${\rm cat}(X,\xi)$, 
{\it (the category of $X$ with respect to class $\xi$)}, 
depending only on the homotopy type of the pair $(X, \xi)$.
It will turn into the classical Lusternik - Schnirelman category ${\rm cat}(X)$ in the case $\xi=0$. 
The main property of $\cat(X,\xi)$ is that it gives a relation between the number of geometrically distinct zeros,
which have closed 1-forms realizing class $\xi$, and the homoclinic orbits of their gradient-like vector fields, see
Theorem \ref{lsmain1}.

Fix a continuous closed 1-form $\omega$ on $X$ representing the cohomology class $\xi$, see the Appendix for
definitions. 

\begin{definition}\label{defcatxi}
{\it We will define ${\rm cat}(X,\xi)$ 
to be the least integer $k$ such that for any integer $N>0$ there exists an open cover
\begin{eqnarray}
X=F\cup F_1\cup\dots \cup F_k,\label{cover}
\end{eqnarray}
such that:
\begin{enumerate}
\item[(a)] Each inclusion $F_j\to X$ is null-homotopic, where $j=1, \dots, k$.
\item[(b)] There exists a homotopy $h_t: F\to X$, where $t\in [0,1]$, such that $h_0$ is the inclusion $F\to X$
and for any point $x\in F$,
\begin{eqnarray}
\int\limits_{\gamma_x}\omega \, \le \, -N,\label{lessthan}
\end{eqnarray}
where the curve $\gamma_x:[0,1]\to X$ is given by $\gamma_x(t) =h_t(x)$.
\end{enumerate}}
\end{definition}

 The meaning of the line integral $\int_\gamma \omega$, where $\gamma:[0,1] \to X$ is
 a continuous curve, is explained in the Appendix.

Intuitively, condition (b) means that in the process of the homotopy $h_t$ 
every point of $F$ makes at least $N$ full twists (in the negative direction)
with respect to $\omega$. We want to emphasize that when $N$ tends to infinity we will obtain
a sequence of different coverings (\ref{cover}) with the same number $k$ and with 
the set $F$ becoming more and more complicated,
so that its limit is wild, looking like a fractal. This explains the approximative nature of our Definition \ref{defcatxi},
which allows to avoid these difficulties.

Observe that $\cat(X,\xi)$ does not depend on the choice of the continuous closed 
1-form $\omega$ (which appears in Definition \ref{defcatxi})
and depends only on the cohomology class $\xi=[\omega]$. 
Indeed, is $\omega'$ is another continuous closed 1-form representing
$\xi$ then $\omega-\omega'=df$, where $f:X\to \R$ is a continuous function, and for any continuous 
curve $\gamma:[0,1]\to X$, 
\[|\int_\gamma \omega - \int_\gamma \omega'| = |f(\gamma(1)) - f(\gamma(0))| \le C,\]
where the constant $C$ is independent of $\gamma$. Here we have used compactness of $X$. 
This shows that if we may construct open covers (\ref{cover}) of $X$
such that (b) is satisfied with an arbitrary large $N>0$ then the same is true with $\omega'$ replacing $\omega$. 

In general the following inequality holds
\begin{eqnarray}
\cat(X,\xi) \le \cat(X),\label{usual}
\end{eqnarray}
since we may always consider covers (\ref{cover}) with $F$ empty. Here $\cat(X)$ denotes the
classical Lusternik - Schnirelman category of $X$, i.e. the least integer $k$ such that
there is an open cover
$
X=F\cup F_1\cup\dots \cup F_k,
$
such that
each inclusion $F_j\to X$ is null-homotopic, where $j=1, \dots, k$.

Inequality (\ref{usual}) can be improved \begin{eqnarray}
\cat(X,\xi) \le \cat(X)-1\label{usual1}
\end{eqnarray}
assuming that $X$ is connected and $\xi\neq 0$. 
This follows because under the above assumptions any open subset $F\subset X$, which is contractible to a point
in $X$, satisfies (b) of Definition \ref{defcatxi} for any $N$. Hence, given a categorical open cover 
$X=G_1\cup \dots\cup G_r$ we may set
$F=G_1$ and $F_j =G_{j+1}$ for $j=1, \dots, r-1$, which gives a cover of $X$ satisfying Definition \ref{defcatxi}.

Observe also that 
\begin{eqnarray}
\cat(X, \xi) \, =\, \cat(X, \lambda \xi), \label{lslambda}
\end{eqnarray}
for $\lambda\in \R$, $\lambda>0$, as follows clearly from the above definition.

\subsection{A reformulation}\label{reformulation}

Sometimes it will be convenient to use a different version of condition (b) of Definition \ref{defcatxi}, 
which we will now describe.

Let $p:\tilde X\to X$ be the normal covering corresponding to the kernel of the homomorphism of periods (cf. Appendix)
\begin{eqnarray}
\pi_1(X,x_0) \to \R, \quad [\gamma]\mapsto \int_\gamma \omega.\label{lsperiods}
\end{eqnarray}
The group of covering transformations $\Gamma$ of this covering equals the image of the homomorphism (\ref{lsperiods}). 
The induced closed 1-form $p^\ast \omega$ equals $df$, where $f: \tilde X\to \R$ is a continuous function. 
Denote $p^{-1}(F)$ by $\tilde F$; it is an open subset invariant under $\Gamma$. Then condition (b) is equivalent to: 
{\it \begin{enumerate}
\item[${\rm (b')}$] There exists a homotopy $\tilde h_t: \tilde F\to\tilde X$, where $t\in [0,1]$,
such that $\tilde h_0$ is the inclusion $\tilde F \to \tilde X$, each $\tilde h_t$ is $\Gamma$-equivariant,
and for any point $x\in \tilde F$,
\[f(\tilde h_1(x)) - f(x) \le - N.\]
\end{enumerate}}
Condition $({\rm b'})$ requires that under the homotopy $h_t$ any point of $\tilde F$ descends at least $N$ units down,
measured by function $f$. The homotopy $\tilde h_t$ is a lift of the homotopy $h_t$ which exists because of the HLP 
(homotopy lifting property) of coverings. If $\tilde \gamma_y$, where $y\in \tilde F$, $p(y)=x$
denotes the path $\tilde \gamma_y(t) =\tilde h_t(y)$ in $\tilde X$
then $\gamma_x = p_\ast\tilde \gamma_y$ and 
\[\int\limits_{\gamma_x} \omega\, =\, \int\limits_{p_\ast \tilde \gamma_y} \omega \, =\,  \int\limits_{\tilde\gamma_y} p^\ast \omega 
\, =
\, \int\limits_{\tilde \gamma_y} df \, =\, f(\tilde h_1(y)) -f(y).\]
This explains equivalence between (b) and $({\rm b'})$.

\subsection{Examples}
\begin{example}\label{lsexample1}
{\rm Consider first the case when $\xi=0$. Let us show that then (\ref{usual}) is an equality.
$\xi=0$ implies $\omega=df$, where $f: X\to \R$ is continuous. Then for any curve $\gamma:[0,1]\to X$ the integral
$\int_\gamma \omega=f(\gamma(1)) - f(\gamma(0))$ 
cannot become smaller than the variation of $f$ on $X$. Therefore, for $\xi=0$ inequality (\ref{lessthan})
may be satisfied for large $N$ only if $F=\emptyset$. This proves $\cat(X,0)=\cat(X)$. }
\end{example}

\begin{example} {\rm
Let $X$ be a mapping torus, i.e. $X$ is obtained from a cylinder
$Y\times [0,1]$, where $Y$ is a compact complex, by identifying points
$(y,0)$ and $(\phi(y),1)$ for all $y\in Y$, where $\phi:Y\to Y$ is a continuous map. Note that we do not assume 
that $\phi$ is a homeomorphism or a homotopy equivalence. 
We will denote points of $X$ by pairs $\langle y,s\rangle$, where $y\in Y$, 
$s\in [0,1]$, understanding that $\langle y,0\rangle =\langle \phi(y),1\rangle $. 
$X$ admits a natural projection $q: X\to S^1$, where $q\langle y, s\rangle = \exp(2\pi i s)$, 
and we will denote by $\omega=q^\ast(d\theta)$ the pullback of the standard angular form $d\theta$ of $S^1$; 
$\omega$ is a closed 1-from on $X$.

Let us show that $\cat(X,\xi)=0$, where $\xi=[\omega]\in H^1(X;\R)$. Given a number $N>0$,
define a homotopy $h_t: X\to X$, where for $t\in [0,1]$
\[h_t\langle y, s\rangle \, =\, \langle \phi^{n(t,s)}(y), s - Nt+n(t,s)\rangle .\]
Here $n(t,s)$ denotes the number of integers contained in the semi-open interval $(s-Nt,s]$. 
It is easy to see that this formula defines a continuous homotopy, $h_0=\id$ and $\int_{\gamma_x}\omega =-N$ 
for any point $x\in X$, where $\gamma_x(t)=h_t(x)$. 
This shows that $\cat(X,\xi)=0$. } 
\end{example}

Note that the Lusternik - Schnirelman category $\cat(X)$ of a mapping torus may be arbitrary large (for example, $\cat(T^n)=n+1$).
Hence the above example shows that the difference 
$$\cat(X)-\cat(X,\xi)$$ 
may be arbitrarily large.

\begin{example} {\rm
In Definition \ref{defcatxi}
the notions \lq\lq up\rq\rq and \lq\lq down\rq\rq appear non symmetrically.
Hence it may happen that $\cat (X,\xi) \neq \cat(X, -\xi)$. We will see such example now.

Consider the mapping torus $X$, as described in the previous example, 
with $Y$ being the sphere $S^2$ and with $\phi: Y\to Y$ a map of degree 2. 
We have seen that $\cat(X,\xi)=0$, where $\xi\in H^1(X;\R)$ denotes the cohomology class described above.
Let us show that $\cat(X, -\xi)\ge 1$.

The universal covering $p: \tilde X\to X$ is Milnor's telescope; it can be described as follows.
$\tilde X$ is obtained from the disjoint union $\amalg_{n\in \Z} X_n$, where $X_n$ denotes $Y\times [n\times n+1]$,
by identifying any point $(y,n)\in X_n$ with $(\phi(y), n)\in X_{n-1}$. The projection $p: \tilde X\to X$ maps any point
$(y,t)\in X_n$ to $\langle \phi^{[t-n]}(y), \{t\}\rangle$. If $\omega$ is the closed 1-form on $X$ described in the previous example
then $p^\ast \omega =df$ where $f:\tilde X\to \R$ is the continuous function given by $f(y,t)=t$.

\setlength{\unitlength}{0.8cm}
\begin{figure}[h]
\begin{center}
\begin{picture}(11,6)
\linethickness{0.1mm}
\put(1,0.5){\vector(1,0){9}}
\put(2,3){\line(0,1){1}}
\put(3.5,2.5){\line(-3,1){1.5}}
\put(3.5,2.5){\line(0,1){2}}
\put(3.5,4.5){\line(-3,-1){1.5}}

\put(3.5,3){\line(0,1){1}}
\put(5,2.5){\line(-3,1){1.5}}
\put(5,2.5){\line(0,1){2}}
\put(5,4.5){\line(-3,-1){1.5}}

\put(5,3){\line(0,1){1}}
\put(6.5,2.5){\line(-3,1){1.5}}
\put(6.5,2.5){\line(0,1){2}}
\put(6.5,4.5){\line(-3,-1){1.5}}

\put(6.5,3){\line(0,1){1}}
\put(8,2.5){\line(-3,1){1.5}}
\put(8,2.5){\line(0,1){2}}
\put(8,4.5){\line(-3,-1){1.5}}

\put(8,3){\line(0,1){1}}
\put(9.5,2.5){\line(-3,1){1.5}}
\put(9.5,2.5){\line(0,1){2}}
\put(9.5,4.5){\line(-3,-1){1.5}}

\put(0.5,3.5){$\tilde X$}
\put(10.5,0.5){$\R$}
\put(5.5,1.25){$\downarrow\quad f$}

\end{picture}
\end{center}
\caption{Covering of mapping torus}
\end{figure}
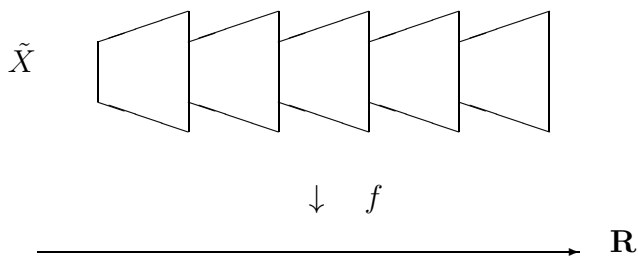

Let $X^k\subset \tilde X$ denote $\bigcup_{n\geq k} X_n$. It is easy to see that $H_2(X^k;\Z)$ is isomorphic to $\Z$ with
the fundamental class of the sphere $Y\times (k+1)\subset X_{k}$ as a generator. The inclusion $X^k\to X^{k-1}$ 
induces on the homology the homomorphism $H_2(X^k;\Z)\to H_2(X^{k-1};\Z)$ of multiplication by $2$. Hence we see that
$H_2(\tilde X;\Z)$ is isomorphic to the abelian group $\Z_{(2)}$
of rational numbers with denominators powers of $2$ and the image of $H_2(X^k;\Z)$
can be identified with $2^{k}\Z\subset \Z_{(2)}$.
It follows that the sphere $Y\times k\subset X_k$ cannot be homotoped into $X^k$ by a homotopy in $\tilde X$. 
Using the condition $({\rm b'})$ of subsection \ref{reformulation}
we see that there is no deformation taking the sphere $Y\times k$ up into $X^k$ and 
this proves that $\cat(X, -\xi)>0$. 

On the contrary, there is a deformation of $\tilde X$ taking all points arbitrarily far down (with respect to $f$), as
we have seen in the previous example.
}\end{example}

\begin{example}\label{lsexample14} {\rm
Consider a bouquet $X=Y\vee S^1$ and assume that the class $\xi\in H^1(X,\R)$ satisfies
$\xi|_Y =0$ and $\xi|_{S^1}\neq 0$. Let us show that then
\[{\rm cat}(X,\xi) \, =\, {\rm cat}(Y)-1.\]
Consider an open cover $X = F\cup F_1\cup\dots \cup F_k$ satisfying conditions (a) and (b) of Definition \ref{defcatxi}.
Let $F'$ denote $F\cap Y$ and $F_j'=F_j\cap Y$. 
We want to show that $F'$ is contractible in $Y$. 
This would imply that $Y=F'\cup F'_1\cup\dots \cup F_k'$ is a categorical cover of $Y$, and hence $\cat(X,\xi)+1\le \cat(Y)$.
Let $h_t: F\to X$ be a homotopy
as in (b) with $N$ large enough. 
For any $x\in F'$, the curve
$h_t(x)$ has to pass along the circle $S^1$ many times and hence it arrives to the point 
$q$ of intersection $S^1\cap Y$ for some $t\in [0,1]$. For 
$x\in F'$ denote by $\sigma(x)\in [0,1]$ the minimal number such that $h_{\sigma(x)}(x)=q$. 
The function $\sigma: F'\to \R$ is continuous. Therefore $h'_t(x)=h_{t\sigma(x)}(x)$ defines a homotopy $F'\to Y$
with $h'_0(x)=x$ and $h'_1(x)=q$.

The opposite inequality $\cat(X,\xi)+1\geq \cat(Y)$ is clear since for any categorical open cover $Y=G_0\cup \dots \cup G_r$
we may set $F=S^1\cup G_0$ and $F_j=G_j$, where $j=1, \dots, r$. The set $F$ satisfies (b) of Definition \ref{defcatxi} for
any integer $N>0$.} 
\end{example}

This example shows that the integer $\cat(X,\xi)$ 
may assume arbitrary non-negative values.

\subsection{Homotopy invariance}

\begin{lemma} Let $\phi: X_1\to X_2$ be a homotopy equivalence, $\xi_i\in H^1(X_i;\R)$, where $i=1, 2$,
and $\xi_1=\phi^\ast(\xi_2)\in H^1(X_1;\R)$.
Then 
\begin{eqnarray}
{\rm cat}(X_1,\xi_1)={\rm cat}(X_2, \xi_2)
\end{eqnarray}
\end{lemma}
\begin{proof} Let $\psi:X_2\to X_1$ be a homotopy inverse of $\phi$. 
Choose a closed 1-form $\omega_1$ on $X_1$ in the 
cohomology class $\xi_1$. 
Then $\omega_2=\psi^\ast \omega_1$ is a closed 1-form on $X_2$ lying in cohomology class $\xi_2$. 

Fix a homotopy $r_t: X_1\to X_1$, where $t\in [0,1]$, such that $r_0=\id_{X_1}$ and $r_1= \psi\circ \phi$.
Compactness of $X_1$ implies that there is a constant $C>0$ such that $|\int_{\alpha_x}\omega_1|< C$
for any point $x\in X_1$, where $\alpha_x$ is the track of the point $x$ under homotopy $r_t$, i.e. 
$\alpha_x(t)=r_t(x)$.

Suppose that ${\rm cat}(X_2, \xi_2)\leq k$. Given any $N>0$, there is an open covering 
$X_2=F\cup F_1\cup\dots\cup F_k$, such that $F_1, \dots, F_k$ are contractible in $X_2$ and 
there exists a homotopy $h_t: F\to X_2$, where $t\in [0,1]$, such that $\int_{\gamma_x}\omega_2\leq -N -C$
for any $x\in F$, where $\gamma_x(t)=h_t(x)$. 
Define 
\[G=\phi^{-1}(F), \quad G_j=\phi^{-1}(F_j), \quad j=1, \dots, k.\]
These sets form an open cover of $X_1=G\cup G_1\cup\dots \cup G_k$. 
Let us show that the set $G\subset X_1$ satisfies condition (b) of Definition
\ref{defcatxi}. Define a homotopy $h'_t: G\to X_1$, where $t\in [0,1]$, by
\begin{eqnarray}
h'_t(x) = \left\{
\begin{array}{ll}
r_{2t}(x), &\mbox{for}\quad 0\le t\le 1/2,\\ \\

\psi(h_{2t-1}(\phi(x))), &\mbox{for}\quad 1/2\leq t \leq 1.
\end{array}
\right.
\end{eqnarray} 
Then $h'_0$ is the inclusion $G\to X_1$ and for any point $x\in G$ holds 
$$\int_{\gamma'_x}\omega_1 \leq -N,$$ 
where $\gamma'_x(t)=h'_t(x)$ is the track of $x$ under homotopy $h_t'$. 

The following diagram 
\begin{eqnarray*}
\begin{array}{ccc}
G_j & \stackrel {\subset}\to & X_1\\
\phi \downarrow && \uparrow \psi\\
F_j & \stackrel{\subset}\to & X_2
\end{array}
\end{eqnarray*}
is homotopy commutative and the horizontal map below is null-homotopic. This shows that the inclusion $G_j\to X_1$
is null-homotopic, where $j=1, \dots, k$. 

The above argument proves that $\cat(X_1,\xi_1) \leq \cat(X_2, \xi_2)$. 

The inverse inequality follows similarly. $\Box$.
\end{proof}

\section{An estimate on the number of zeros}

In this section we will use invariant $\cat(X,\xi)$ to obtain a lower bound on the number of
zeros of vector fields having no homoclinic orbits.

Let $\omega$ be a smooth closed 1-form on a connected closed smooth manifold $M$. 
We will assume that $\omega$ has finitely many zeros.
A smooth vector field $v$  is a {\it gradient-like} vector field for 
on $M$ if
(1)
$\omega(v)>0$ 
on the complement of the set of zeros of $v$, and 
(2) in a neighborhood $U_p\subset M$
of any zero $p$, the field $v|_{U_p}$ coincides with the gradient of 1-form $\omega$ 
with respect to a Riemannian metric on $U_p$. 

An integral trajectory $\gamma(t)$ of vector field $v$ is {\it a homoclinic orbit} if both limits 
$\lim_{t\to \pm \infty} \gamma(t)=\gamma(\pm \infty)$
exist and are equal. The point $p=\gamma(\pm \infty)$ is then a zero of $v$. 
Note that for
 any homoclinic orbit $\gamma$ holds 
$\int_\gamma \omega >0$. 
Hence, homoclinic orbits do not exist in the case $\xi=0$, i.e. in the usual gradient systems.

More generally, a {\it homoclinic cycle of length $n$} is a sequence of trajectories $\gamma_1, \dots, \gamma_n$
of the field $v$ such that all the limits $\lim_{t\to \pm \infty} \gamma_i(t)$ exist, where $i=1, \dots, n$, and 
\[\lim_{t\to + \infty} \gamma_i(t) = \lim_{t\to - \infty} \gamma_{i+1}(t), \quad i=1, \dots, n.\]
For $i=n$ this means that $\lim_{t\to + \infty} \gamma_n(t) = \lim_{t\to - \infty} \gamma_1(t)$, and so the union of the 
curves $\gamma_i$ form a closed cycle.

\begin{theorem}\label{lsmain1} Let $\omega$ be a smooth closed 1-form on a closed manifold $M$ and let 
$\xi=[\omega]\in H^1(M;\R)$ denote the cohomology class of $\omega$.
If $\omega$ admits a gradient-like vector field $v$ with no homoclinic cycles, then $\omega$ has at least 
$\cat(M,\xi)$ geometrically distinct zeros.
\end{theorem}

Here is a different formulation of the above Theorem:

\begin{theorem}\label{lsmain2} 
If the number of zeros of a smooth closed 1-form $\omega$ is less than $\cat(M,\xi)$, where
$\xi=[\omega]\in H^1(M;\R)$ denotes the cohomology class of $\omega$,
 then any
 gradient-like vector field for $\omega$ has a homoclinic cycle.
\end{theorem}

Combined with Theorem \ref{lsreal}, this shows that there may 
exist homoclinic cycles which cannot be destroyed while perturbing the gradient-like vector field.
This {\it \lq\lq focusing effect\rq\rq} starts when the number
of zeros of a closed 1-form becomes less than $\cat(M,\xi)$. 
It is a new phenomenon, not occuring in the Novikov theory: assuming that 
the zeros of $\omega$ are all Morse
it is always possible 
to find a gradient-like vector fields $v$ for $\omega$
such that any integral trajectory connects a zero with higher Morse index with a zero with lower Morse index
(by the Kupka - Smale Theorem \cite{Sm}).

\subsection{Proof of Theorem \ref{lsmain1}}

Let $p_1, \dots, p_k\in M$ denote all the zeros of $\omega$. 
Assume that there exists a gradient-like vector field $v$ for $\omega$ with no homoclinic orbits.
Our purpose is to show that then $\cat(M,\xi)\leq k$.

Fix a number $N>0$.

Consider the flow $M\times \R\to \R$, where $(m,t)\mapsto m\cdot t$, generated by the field $-v$. 

Choose small closed disks 
$U_j$ around each point $p_j$, where $j=1, \dots, k$. 
We will assume that for $i\neq j$ the disks $U_i$ and $U_j$ are disjoint. Also, we will fix a Riemannian metric on $M$
such that in the disks $U_i$, where $i=1, \dots, k$,
the field $v$ is the gradient of $\omega$ with respect to this metric.

We claim that: {\it one may 
choose closed disks $V_i$, where $i=1, \dots, k$, 
such that: 

(a) $p_i\in {\rm Int}\, V_i$ and $V_i\subset {\rm Int}\, U_i$, 

(b) the disk $V_i$ is gradient-convex (see Appendix \ref{lsapp2}) in the following sense.
Consider the covering $\pi: \tilde M\to M$ corresponding to the kernel of the homomorphism of periods
(\ref{lsperiods}). Then the form $\omega$ lifts to $\tilde M$ as a smooth function $f: \tilde M\to M$,
i.e. $\pi^\ast \omega =df$, and the field $v$ lifts to a 
gradient-like vector field $\tilde v$ of function $-f$. We require any lift to $\tilde M$ of the disk $V_i$ to be 
gradient-convex, see Appendix \ref{lsapp2}.

(c) Let $\partial_- V_i$ denote the set of points $p\in \partial V_i$, such that for all sufficiently small 
$\tau >0$, holds $p\cdot \tau \notin V_i$. 
Then we require that for any $p\in \partial_-V_i$ there exists no real number $t_p>0$, 
such that: the point $p\cdot t_p$ belongs to ${\rm Int}\, V_i$, and 
$\int_{p}^{p\cdot t_p}\omega \, \geq \, -N,$
where the integral is calculated along the integral trajectory $\sigma_p: [0,t_p]\to M$,
$\sigma_p(t)=p\cdot t$. }

Assume first that $V_i$ is any neighborhood satisfying (a) and (b). Then any trajectory 
$p\cdot t$, where $p\in \partial_-V_i$, leaves $U_i$ before it may re-enter $V_i$. This follows from gradient
convexity of $V_i$, since the disk 
$U_i$ also lifts to the covering $\tilde M\to M$. Hence:

(i) {\it There exists $a>0$, such that for any $p\in \partial_-V_i$ and $t>0$ with $p\cdot t\in V_j$ holds $t\geq a$.}

We may take $a=\min\{l_i v_i^{-1}; i=1, \dots, k\}$, 
where $l_i>0$ denotes the distance between $V_i$ and $M-{\rm Int}\, U_i$, and 
$v_i =\max |v(x)|$ for $x\in U_i -{\rm Int}\, V_i$. 

Note that if we shrink the disks $V_i$ the number $a>0$ may only increase, assuming that $V_i$ are sufficiently small.

(ii) {\it There exists $b>0$, such that for any $p\in \partial_-V_i$ and $t>0$ with $p\cdot t\in V_j$ holds} 
$$\int_p^{p\cdot t}\omega < -b.$$ 

We may take $b=\min\{\epsilon_i a_i; i=1, \dots, k\}$, where $-\epsilon_i=\max \omega(v)(x)$ for $x\in U_i-{\rm Int}\, V_i$.

Suppose that we may never achieve (c) by shrinking the disk
$V_i$ satisfying conditions (a) and (b). 
Then there exists a sequence of points $p_{i,n}\in \partial_- V_i$ 
and two sequences
of real numbers $t_{i,n}>0$ and $s_{i,n}<0$ such that:

(1) the set $p_{i,n}\cdot [s_{i,n},0]$ is contained in the disk $V_i$ and 
the point $p_{i,n}\cdot s_{i,n}$ converges to $p_i$ as $n$ tends to $\infty$;

(2) $p_{i,n}\cdot t_{i,n}$ converges to $p_i$;

(3) $$\int\limits_{p_{i,n}}^{p_{i,n}\cdot t_{i,n}}\omega \geq -N.$$

Let $\gamma_{i, n}\in H_1(M)$ denote the homology class of the loop obtained as follows.
Start at the point $p_{i,n}\in \partial_-V_i$, follow the trajectory $p_{i,n}\cdot [0,t_{i,n}]$, and then connect the endpoint
$p_{i,n}\cdot t_{i,n}$ with $p_{i,n}$ by a path inside $V_i$. 
We claim that {\it the set $\{\gamma_{i,n}\}_{n\geq 1}\subset H_1(M)$ of 
the obtained homology classes is finite.} This would follow once we show that the total length $L_{i,n}$
of the parts of any trajectory
$p_{i,n}\cdot [0,t_{i,n}]$ spent outside the set $U_1\cup \dots\cup U_k$, is bounded. 
Writing
\[\int_\gamma \omega \, =\, \int_a^b \omega(\dot\gamma(t))dt \, =\, \int_a^b \frac{\omega(\dot\gamma(t))}{|\dot\gamma(t)|}\cdot
|\dot\gamma(t)|dt\]
and using (3) above, it is easy to see that 
\[L_{i,n} \leq N\cdot c^{-1},\] 
where $c>0$ is given by
\[ c\, =\, - \max\, \{\omega(v_x)\cdot |v_x|^{-1}\, ;\,  x\in M - \bigcup_{i=1}^k {\rm Int}\, U_i\}.\]

Passing to a subsequence, we may assume that $p_{i,n}$ converges
to a point $q_i\in \partial_- V_i$ and 
the sequences $s_{i,n}$ and $t_{i,n}$ have finite or
infinite limits $s_{i}$ and $t_{i}$ correspondingly. 
We may also assume that the homology class $\gamma_{i,n}\in H_1(M)$ is
independent of $n$. 

 (1) and (2) imply $t_{i}=+\infty$ and $s_i =-\infty$. 

First we want to show that 
$\lim_{t\to - \infty} q_i\cdot t =p_i.$ 
We will identify $V_i$ with one of its lifts to the covering $\tilde M$. 
If there exists $s<0$ such that $f(q_i\cdot s)>f(p_i)+\epsilon$ for some $\epsilon >0$, then the set 
$\{(x,s)\in \tilde M\times \R; f(x\cdot s)>f(p_i)+\epsilon\}$ (which is open) contains $(p_{i,n},s_{i,n})$ for all
$n$ sufficiently large, which contradicts (1). This shows that on the covering $\tilde M$ holds
$q_i\cdot s\to p_i$ for $s\to -\infty$. Hence the same relation holds on the initial manifold $M$, as well. 

Now we want to understand the limit $\lim_{t\to +\infty}q_i\cdot t$.
Consider the above lift of $V_i$ to the covering $\tilde M$, such that the flow becomes the gradient-like flow 
$\tilde v$ of function
$-f:\tilde M\to \R$.
The points $p_{i,n}\cdot t_{i.\,n}$ all belong to a translate $g V_i$ of the disk $V_i$, where $g$ is independent of $n$. 
The trajectory $q_i\cdot t$ in $\tilde M$ for large $t$ may either reach the neighborhood $gV_i$, or it may
be "caught" by some other critical point of $f$ on the way.

Let us show that in the first case the point $q_i\cdot t$ tends to $g p_i$, as $t$ tends to $+\infty$ and hence the vector field
$v$ on $M$ has a homoclinic orbit, starting and ending at $p_i$.  
If it is not true, then $f(q_i\cdot t)< f(g p_i)-\epsilon$ for some $t$ and $\epsilon >0$.
Then the
open set 
$\{(x,t)\in \tilde M\times \R; f(x\cdot s)<f(g p_i)-\epsilon\}$ would contain $(p_{i,n},t_{i,n})$ for all
$n$ sufficiently large, which contradicts  (2). This shows that on the covering $\tilde M$ holds
$q_i\cdot t\to g p_i$ for $t\to +\infty$. Hence on $M$, the trajectory $q_i\cdot t$ tends to $p_i$ for $t\to +\infty$.  

Consider now the second possibility; we will show that in this case $v$ has a homoclinic cycle of length greater than 1.
Assume that the limit $q_i\cdot t$ for $t\to \infty$ equals $hp_j$, where $h$ is an element of the group of periods of $\omega$,
and 
$$f(gp_i)+b < f(hp_j) \quad \mbox{and}\quad  f(hp_j) < f(p_i)-b.$$ 
Then for large $n$, the trajectory starting at $p_{i,n}$ enters the neighborhood $hV_{j}$
and leaves it at some point $p'_{i,n}\in \partial_- (hV_j)$. Passing to a subsequence, we may assume that the
sequence $p'_{i,n}$ converges to a point $q'_i\in \partial_-(hV_j)$. Then, as above, $\lim_{t\to -\infty} q'_i\cdot t$
equals $hp_j$ and the limit $\lim_{t\to +\infty} q'_i\cdot t$ either equals $gp_i$, or it equals 
$h_1p_{j_1}\in \tilde M$, where
$$f(gp_i)+b< f(h_1p_{j_1})\quad\mbox{and}\quad f(h_1p_{j_1})<  f(hp_j)-b.$$ 
Continuing these arguments by induction, we obtain a homoclinic cycle on the manifold
$M$, \lq\lq starting and ending\rq\rq at $p_i$. 
The number of steps in the above process is finite (at most $[N/b]$). 

The above argument proves existence of the disks $V_i$ with properties (a), (b), (c), assuming that vector field $v$
has no homoclinic cycles.

Next we will construct an open cover $F\cup F_1\cup F_2\cup \dots \cup F_k =M$. 

We will define $F$ as the set of all points $p\in M$ such that there exists a positive number $t_p>0$, 
so that the integral curve $\sigma_p: [0,t_p]\to M$, where
$\sigma_p(t)=p\cdot t$, satisfies $\int_{\sigma_p}\omega =-N$.
It is clear that $F$ is open, and $p\mapsto t_p$ is a continuous real valued function on $F$. 

We may define a homotopy
$$h_\tau: F \to F, \quad\mbox{by}\quad h_\tau(p) =p\cdot (\tau t_p),\quad \tau\in [0,1]$$ 
This homotopy satisfies condition (b) of Definition
\ref{defcatxi}.

Now we will define the sets $F_j$, where $j=1, \dots, k$. 
We say that $p\in F_j$ if for some $t_p>0$ the point $p\cdot t_p$ belongs to the interior of $V_j$, 
and
$$\int\limits_{\sigma_p}\omega >-N,$$
 where
$\sigma_p: [0,t_p]\to M$ is given by $\sigma_p(t)=p\cdot t$. 
It is clear that $F_j$ is open. 

The sets $F, F_1, \dots, F_k$ cover $M$. Indeed, 
for any point $p\in M$ either 
$$\int\limits_p^{p\cdot t} \omega < -N$$
for some
$t>0$, or 
the trajectory 
$\gamma(t)=p\cdot t$ \lq\lq enters\rq\rq\,  a zero $p_j$, so that
$$\lim_{t\to \infty}\, \int\limits_p^{p\cdot t} \omega \geq  -N.$$
In the first case $p$ belongs to $F$, and in the second case $p$ belongs to $F_j$.

Now we will show that the set 
$F_j$, where $ j=1, \dots, k$,
 is contractible in $M$. For any point $p\in M$ let $J_p\subset \R$ denotes the set
$J_p =\{t\geq 0; p\cdot t\in V_j\}$.
Because of our assumption about gradient-convexity of $V_j$, the set $J_p$ is a union of disjoint closed intervals,
and some of these intervals may degenerate to a point. 
Consider the first interval $[a_p,b_p]\subset J_p$. 
If this interval degenerates to a point (i.e. the trajectory through $p$ touches $V_j$), 
then $p$ does not belong to the set $F_j$, according to our assumption (c). By the same reason, points of 
$\partial_-V_j$ do not belong to $F_j$. 

Assume now $p\in F_j$ and $p\notin {\rm Int}\, V_j$.
Then the point $p\cdot t$ lies in the interior of $V_j$ for $a_p<t<b_p$. Also, we have 
$$\int\limits_p^{p\cdot a_p}\omega > -N.$$
The function 
$\phi_j: F_j\to \R$, given by
$$
\phi_j(p) = \left\{
\begin{array}{ccl}
0,&\mbox{for}& p\in {\rm Int}\,  V_j,\\ 
a_p,&\mbox{for}&p\in F_j -{\rm Int}\, V_j
\end{array}
\right.
$$
is continuous. To show this, suppose that a sequence of points $x_n\in F_j$, where $n=1, 2, \dots$,
 converges to $x_0\in F_j$.
First consider the case $x_0\in V_j$. 
Since $x_0\in F_j$, we conclude that $x_0\notin \partial_-V_j$ (because of condition (c)), i.e.
$x_0\cdot \tau$ belongs to ${\rm Int}\,V_j$ for all $\tau \in (0, \epsilon)$. It follows that 
for any $\tau \in (0, \epsilon)$ holds $x_n\cdot \tau \in {\rm Int}\, V_j$ for all large $n$. Therefore,
$a_{x_n}<\tau$ for all large $n$. Hence, the sequence $a_{x_n}$ converges to $0$. 
Consider now the case when $x_0$ does not belong to $V_j$. Then the trajectory $x_0\cdot t$ does not touch $V_j$
for $t<a_{x_0}$ (again, because of condition (c))
and $t=a_{x_0}$ is the first moment when the trajectory penetrates $V_j$. We know also that
the velocity vector $v_x$ is transversal to the boundary $\partial V_j$, where $x=x_0\cdot a_{x_0}$. 
Then the sequence $a_{x_n}$ converges to $a_{x_0}$,
 as follows from continuity of solutions of ordinary differential equation with respect to the initial conditions.

We may define a homotopy 
$$h_\tau: F_j \to M, \quad h_\tau(p) = p\cdot (\tau\phi_j(p)), \quad p\in F_j, \quad \tau \in [0,1].$$
Here $h_0$ is the inclusion $F_j\to M$ and $h_1$ maps $F_j$ into the disk $V_j$. 
 
This completes the proof. $\Box$

\section{Moving homology classes}

Here we study the effect of condition ${\rm (b')}$ of subsection \ref{reformulation}
on homology classes. For simplicity we assume that the group of periods $\Gamma\subset \R$ is infinite cyclic.

$\kk$ will denote a field. 

Let $\tilde X\to X$ be an infinite cyclic covering, i.e. a regular covering of a finite CW-complex having an
infinite cyclic group of covering transformations. In this section 
we will denote by $\tau:\tilde X \to \tilde X$ a fixed generator of this group.

Let $K\subset \tilde X$ be a compact subset such that $\tilde X$ is the union of the translates
$\tau^i (K)$, where $i\in \Z$.

\begin{definition}{\it We will say that a homology class $z\in H_q(\tilde X;\kk)$ is movable to $+\infty$
if for any integer $N\in \Z$ there exists a cycle in $\bigcup_{i>N}\tau^i (K)$ representing $z$.
Similarly, a homology class $z\in H_q(\tilde X;\kk)$ is movable to $-\infty$
if for any integer $N\in \Z$ there exists a cycle in $\bigcup_{i<N}\tau^i (K)$ representing $z$.}
\end{definition}

It is clear that the above properties of homology class $z\in H_q(\tilde X;\kk)$
do not depend on the choice of compact $K$.

The following Lemma gives an approximative condition of movability.
Roughly, it claims: if a cycle can be moved sufficiently large distance away then it may be moved 
arbitrarily far away. 
The word \lq\lq {\it cycle}\rq\rq means \lq\lq {\it singular cycle with coefficients in $\kk$}\rq\rq.

\begin{lemma}\label{lsmov} Let $K\subset \tilde X$ be a compact subset such that 
$\bigcup_{j\in \Z} \tau^j (K)$ coincides with $\tilde X$. 
Then there exists an integer $N>0$ (depending on $K$),
such that the following properties (a) and (b) hold:

(a) Let a cycle $c$ in $K$ be homologous in $\tilde X$ 
to a cycle in $\bigcup_{j \geq N} \tau^j(K)$. Then the 
 homology class $[c]\in H_q(\tilde X;\kk)$ is 
movable to $+\infty$; 

(b) Let a cycle $c$ in $K$ be homologous in $\tilde X$ 
to a cycle in $\bigcup_{j \leq -N} \tau^j(K)$. Then the
 homology class $[c]\in H_q(\tilde X;\kk)$ is 
movable to $-\infty$.
\end{lemma} 
\begin{proof} Denote  
\[V_r = \im[ H_q( K;\kk) \to H_q(\tilde X;\kk)]\cap \im[ H_q( \left(\bigcup_{j\geq r}\tau^{j} K\right);\kk) \to H_q(\tilde X;\kk)].\]
where $r=1, 2, \dots$.
Then $V_1\supset V_2\supset V_3\supset \dots$ is a decreasing sequence of finite-dimensional vector spaces.
Hence, there exists
an integer $N$, such that $V_N$ coincides with $V_\infty = \cap_{r>0} V_r$. 
Any homology class $z\in V_\infty$ is movable to $+\infty$. Thus, this $N$ satisfies (a). 
By the similar reasons, we may increase $N$, if necessary, such that (b) is satisfied as well. $\Box$
\end{proof}

\begin{lemma}\label{lslm31} Given a homology class $z\in H_q(\tilde X;\kk)$, the following conditions are equivalent:
\begin{enumerate}
\item[(i)] $z$ is movable to $+\infty$;
\item[(ii)] $z$ is movable to $-\infty$;
\item[(iii)] $z$ is a torsion element of $\kk[\tau, \tau^{-1}]$-module $H_q(\tilde X;\kk)$, i.e. there exists a nontrivial 
Laurent polynomial
$p(\tau)\in \kk[\tau, \tau^{-1}]$, such that $p(\tau)z=0$.
\end{enumerate}
\end{lemma}

\begin{proof} We will show that (i) implies (iii) and that (iii) implies both (i) and (ii). 
The implication (ii) $\Rightarrow$ (iii) follows similarly.

Assume that a class $z\in H_q(\tilde X;\kk)$ is movable to $+\infty$. Realize $z$ by a cycle $c$ in $\tilde X$
and specify a compact subset $K\subset \tilde X$ containing $c$ and such that $\tilde X = \bigcup_{i\in \Z} \tau^i (K)$. 
Assume that $N>0$ is large enough, so that it satisfies Lemma \ref{lslm31}
and the subset $\tau^N (K)$ is disjoint from $K$. 
Let us show that for any integer $r\geq N$
the homology class $z$ may be realized by a cycle in $\tau^r (K)$. 
Since $z$ is movable to $+\infty$, 
there exists a cycle $c'$ in $U=\bigcup_{j\geq 2r}\tau^j(K)$ representing $z$. 
Write $\tilde X = B\cup C$, where $B$ contains $K$, $C$ contains $U$ and $B\cap C=\tau^r K$.
In the Mayer - Vietoris sequence
\[H_q(B\cap C;\kk) \to H_q(B;\kk)\oplus H_q(C;\kk)\to H_q(\tilde X;\kk),\]
the difference $[c] - [c']$ goes to zero. Hence there is a cycle in $B\cap C=\tau^r K$, which is homologous to $c$ in $B$
and homologous to $c'$ in $C$. This proves our claim.

Consider
\[V = \bigcap_{r\ge N} \im  [H_q(\tau^r K;\kk) \to H_q(\tilde X;\kk)]\, \subset H_q(\tilde X;\kk).\]
It is a finite dimensional $\kk$-linear subspace; we observe that $V$ 
is invariant under $\tau^{-1}$. Hence, by the Caley - Hamilton Theorem,
there exists a polynomial $p(\tau)\in \kk[\tau]$, such that $p(\tau^{-1})$ acts trivially on $V$. Since $z$ belongs to $V$,
we obtain $p(\tau^{-1})z=0$. This shows that (i) implies (iii). 

Let us show that (iii) implies (i) and (ii). Suppose that $z\in H_q(\tilde X;\kk)$ is such that $p(\tau)z=0$ for a Laurent polynomial
\[p(\tau) = \sum\limits_{i=r}^{r+\ell}a_i\tau^i, \quad a_i\in \kk, \quad \mbox{where}\quad a_r\neq 0, \quad a_{r+\ell}\neq 0.\]
Consider the ring $R = \kk[\tau, \tau^{-1}]/J$, where $J$ is the ideal generated by $p(\tau)$. The powers $\tau^i$, where
$i=r, r+1, \dots, r+\ell -1$, form an additive basis of $R$. Since multiplication by $\tau$ is an automorphism of $R$,
we obtain that for any integer $N$, 
the powers $\tau^N, \tau^{N+1}, \dots, \tau^{N+\ell-1}$ form a linear basis of $R$ as well. In particular, we
may express $1\in R$ as a linear
combination of $\tau^N, \tau^{N+1}, \dots, \tau^{N+\ell -1}$ in $R$. 
This means that for any $N$ we may find numbers
$b_j\in \kk$, where $j=N, \dots, N+\ell-1$, such that
\[ z = \sum_{j=N}^{N+\ell-1} b_j\tau^j z\quad \mbox{in}\quad H_q(\tilde X;\kk).\]

Assume that $K\subset \tilde X$ is a compact subset such that 
$\tilde X = \bigcup_{j\in \Z}\tau^j K$
and class $z$ can be realized by a cycle $c$ in $K$. Then for any integer
$N$ class $z$ may be realized by the cycle
$$\sum_{j=N}^{N+\ell-1} b_j\tau^j c \quad \mbox{lying in}\quad \bigcup_{N\leq j < N+\ell }\tau^j K.$$ 
Hence $z$ is movable to both ends $\pm\infty$ of $\tilde X$. $\Box$
\end{proof}

\section{Cohomological lower bound for $\cat(X,\xi)$}

In this section we will give cohomological lower bounds on $\cat(X,\xi)$. 
For simplicity we assume that $\xi$ is integral, i.e. $\xi\in H^1(X;\Z)$.

\subsection{Statement of the result}

Let $\kk$ be a field. 

Given a finite CW-complex $X$ and an integral cohomology class $\xi\in H^1(X;\Z)$. For any nonzero $a\in \kk$ 
there is a local system 
over $X$ with fiber $\kk$ such that the monodromy along any loop $\gamma$ in $X$ is 
multiplication by $a^{\langle \xi, \gamma\rangle}$
$$a^{\langle \xi, \gamma\rangle}: \kk \to \kk, \quad q\mapsto a^{\langle \xi, \gamma\rangle}q,
\quad \mbox{for}\quad q\in \kk.$$ 
Here $\langle \xi, \gamma\rangle\in \Z$ denotes the value of the class $\xi$ on the loop $\gamma$.
This local system will be denoted $a^\xi$. The cohomology of this local system $H^q(X;a^\xi)$ 
is a vector space over $\kk$ of finite dimension. Note that for $a=1$ the local system $a^\xi$ is the constant local system 
$\kk$. If $a, b\in \kk^\ast$ are two nonzero numbers, there is an isomorphism of local systems 
$$a^{\xi}\otimes b^{\xi} \simeq (ab)^{\xi}.$$
Hence we have well-defined cup-product pairing
\begin{eqnarray}
\cup: \quad H^q(X;a^\xi)\otimes H^{q'}(X;b^\xi) \to H^{q+q'}(X; (ab)^\xi).
\end{eqnarray}

The following Theorem is one of the main results of this section.

\begin{theorem}\label{lscup} Assume that there exist cohomology classes 
$$u\in H^{q}(X;a^\xi),\quad v\in H^{q'}(X;b^\xi),\quad w_j\in H^{d_j}(X;\kk),\quad\mbox{where}\quad  j=1, \dots, r,$$
such that 

{\rm (i)} $d_1>0, \dots, d_r>0$; 

{\rm (ii)} the cup product
\begin{eqnarray}
u\cup v\cup w_1\cup w_2\cup\dots \cup w_r\in H^d(X;(ab)^\xi), \label{lscup1}
\end{eqnarray}
is nontrivial, where $d =q+q' +d_1 +\dots +d_r$;

{\rm (iii)} the numbers $a, b\in \kk^\ast$ do not belong to a finite subset ${\rm Supp}(X,\xi)\subset \kk^\ast$ depending on 
the pair $(X,\xi)$, see below. 

Then $\cat(X,\xi)>r$.
\end{theorem}

Property (iii) may be also expressed by saying that the numbers $a\in \kk^\ast$ and $b\in \kk^\ast$ are {\it generic}.

For $\xi\neq 0$ and $a\ne 1$ the bundle $a^\xi$ admits no globally defined flat sections and so 
the zero-dimensional cohomology $H^0(X;a^\xi)=0$ vanishes. Hence for $\xi\ne 0$ the classes $u$ and $v$ 
in Theorem \ref{lscup} must automatically have positive degrees, i.e. $q>0$ and $q'>0$.

\subsection{The set ${\rm Supp}{(X,\xi)}$}
For $\xi=0$ we will define ${\rm Supp}(X,\xi)\subset \kk^\ast$ to be the empty set.

Now we will define the set ${\rm Supp}(X,\xi)\subset \kk^\ast$ assuming that $\xi\in H^1(X;\Z)$ is nonzero and indivisible.
Let $p:\tilde X\to X$ be the infinite cyclic covering corresponding to $\xi$. 
Fix a generator $\tau: \tilde X\to \tilde X$ of the group of covering translations such that 
$f(\tau x) =f(x) +1$ for all $x\in \tilde X$, where $f:\tilde X\to \R$ is a continuous function satisfying $df=p^\ast\omega$.
Here $\omega$ is a closed 1-form on $X$ in the cohomology class $\xi$. Since the translate $\tau: \tilde X\to \tilde X$
acts on the homology, we obtain that $H_\ast(\tilde X;\kk)$ is 
a graded module over the ring $\Lambda =\kk[\tau, \tau^{-1}]$ of Laurent polynomials. The ring $\Lambda$ is Noetherian,
hence $H_\ast(\tilde X;\kk)$ is finitely generated over $\Lambda$. Its $\Lambda$-torsion submodule
$$T = {\rm Tor}_{\Lambda}(H_\ast(\tilde X;\kk))$$
(which, according to the previous section, coincides with the set of homology classes movable to $\pm \infty$) 
is a finite-dimensional 
$\kk$-vector space. Multiplication by $\tau$ is an invertible $\kk$-linear operator $\tau: T\to T$. We will define the 
set ${\rm Supp}(X,\xi)\subset \kk^*$ as the set of eigenvalues of $\tau^{-1}: T\to T$. 

In the case when $\xi=\lambda\eta$, where $\eta\in H^1(X;\Z)$ is indivisible and $\lambda\in \Z$, $\lambda>0$, 
we will define the set ${\rm Supp}(X,\xi)$
as $\{a\in \kk^\ast; a^\lambda\in {\rm Supp}{(X,\eta)}\}$. 

Note that for $\kk=\C$ the set ${\rm Supp}(X,\xi)\subset \C$ consists of finitely many algebraic numbers. 
Hence for $\kk=\C$ in Theorem \ref{lscup}
one may always take for $a, b\in \C^\ast$ arbitrary transcendental numbers. 

\subsection{Lifting property}
The proof of Theorem \ref{lscup}
is based on the following lifting property of cohomology classes:

\begin{theorem}\label{lslift}
Let $X$ be a finite CW-complex, $\xi\in H^1(X;\Z)$,  $a\in \kk^\ast$, $a\notin {\rm Supp}(X,\xi)$. 
Fix a continuous closed 1-form $\omega$ in cohomology class
$\xi$ and a compact $K\subset X$, such that $\xi|_K =0$. 
Then for any open subset $F\subset K$ satisfying condition (b) of Definition
\ref{defcatxi} with respect to closed 1-form $\omega$ with integer $N>0$ large enough, the homomorphism
\begin{eqnarray}
H^q(X,F;a^\xi) \to H^q(X;a^\xi)\label{lslifting}\label{lsepi}
\end{eqnarray}
is an epimorphism.
\end{theorem}

{\bf Proof of Theorem \ref{lslift}.} We already observed in Example 
\ref{lsexample1} that for $\xi=0$ condition (b) of Definition \ref{defcatxi}
may be satisfied for large $N$ only if $F=\emptyset$. Therefore
Theorem \ref{lslift} is true for $\xi=0$. 

We will assume below that $\xi\not= 0$. 
Moreover, without loss of generality we may assume that class $\xi\in H^1(X;\Z)$ is indivisible. 

Let $f: \tilde X \to \R$ be a continuous function (unique up to a constant) satisfying $df=p^\ast \omega$. 
Then $f(\tau x)=f(x)+1$ for all $x\in \tilde X$, where $\tau: \tilde X \to \tilde X$ a generator of the group of translations. 
Since $\xi|_K=0$, the compact $K\subset X$ may be lifted to $\tilde X$. Fix such lift $K\subset \tilde X$.
We will assume that the function $f|_K$ assumes values in $[0,c]\subset \R$ for some integer $c>0$;
this may be always achieved by adding a constant to $f$. 

Let $N'>0$ be the number given by Lemma \ref{lsmov}
applied to a lift of $K$ to $\tilde X$. Then any cycle in $K\subset \tilde X$, which is homologous in $\tilde X$
to a cycle in $\bigcup_{j\leq -N'} \tau^j(K)$, is movable to $-\infty$, 
and hence (according to Lemma \ref{lslm31}) 
represents a $\Lambda$-torsion homology class
in $H_\ast(\tilde X;\kk)$.  

It follows, that for any subset $F\subset K$, satisfying condition ${\rm (b')}$ of subsection \ref{reformulation}
with $N>N'$,
the homomorphism $H_\ast(\tilde F;\kk) \to H_\ast(\tilde X;\kk)$ 
induced by the inclusion $\tilde F\to \tilde X$, takes values in $\Lambda$-torsion submodule $T$
of $H_\ast(\tilde X;\kk)$. 
The set $\tilde F$ is a disjoint union of infinitely many copies of $F$ and hence 
the homology of $\tilde F$ is $H_\ast(\tilde F;\kk)\simeq H_\ast(F;\kk)\otimes_{\kk}\Lambda$.

The claim that (\ref{lsepi}) is an epimorphism is equivalent to the claim that $H^\ast(X;a^\xi)\to H^\ast(F;a^\xi)$ 
is the zero map. 
Using duality between homology and cohomology we see the later is equivalent to the statement 
that the homomorphism $H_\ast(F;a^{-\xi}) \to H_\ast(X;a^{-\xi})$,
induced by the inclusion $F\to X$ on the homology of the dual local system $a^{-\xi}$, is zero. 
Since $\xi|_F =0$, the local system $a^{-\xi}|_F\, \simeq \, \kk$ is trivial. Hence we want to show that the 
homomorphism $H_\ast(F;\kk)\to H_\ast(\tilde X;a^{-\xi})$ is zero. 
The inclusion $F\to X$ equals the composition
$F\stackrel {\subset}\to K\to \tilde X\stackrel p\to X$, and we know that $H_\ast(F;\kk)\to 
H_\ast(\tilde X;\kk)$ takes values in the $\Lambda$-torsion submodule $T$. 

To complete the proof it is enough to show that for $a\notin {\rm Supp}(X,\xi)$ holds
$p_\ast(T)=0$, where the homomorphism 
$p_\ast: H_\ast(\tilde X;\kk)\to H_\ast(X;a^{-\xi})$ is induced by the covering projection $p: \tilde X\to X$. 
Consider the following well-known exact sequence
\begin{eqnarray*}
 \dots \to H_i(\tilde X;\kk) \stackrel{\tau-b}\longrightarrow H_i(\tilde X;\kk) \stackrel{p_\ast}\to H_i(X;b^{\xi})\to \dots,
\end{eqnarray*}
where $b=a^{-1}$.
The linear map $\tau -b: T\to T$ is an isomorphism and hence the submodule 
$T\subset H_\ast(\tilde X;\kk)$ is contained
in the image of $\tau -b$ and therefore $p_\ast(T)=0$. 
This completes the proof. 
$\Box$

\subsection{Proof of Theorem \ref{lscup}}
If $\xi=0$ under the conditions of Theorem \ref{lscup} the classical Lusternik - Schnirelman theory
gives $\cat(X) >r+2$; hence Theorem \ref{lscup} holds for $\xi=0$.

We will assume below that the class $\xi\not=0$ is nonzero and indivisible.

Suppose that $\cat(X,\xi)\le r$. Let us show that any cup-product (\ref{lscup1}) satisfying conditions
(i), (ii), (iii) of Theorem \ref{lscup} must vanish. Since $\xi$ is an integral class $\xi\in H^1(X;\Z)$, 
we may cover $X$ by two open subsets 
$$X=U\cup W,\quad\mbox{such that}\quad \xi|_{\bar U} =0, \quad \xi|_{\bar W}=0.$$
Choose a continuous closed 1-form $\omega$ in cohomology class $\xi$. 
Our assumption $\cat(X,\xi)\le r$ implies that for any $N>0$ there is an open cover 
$X=F\cup F_1\cup \dots \cup F_r$ satisfying properties 
(a) and (b) of Definition \ref{defcatxi}.

Each cohomology class $w_j\in H^{d_j}(X;\kk)$ 
may be lifted to a relative cohomology class $\tilde w_j\in H^{d_j}(X,F_j;\kk)$, where $j=1, \dots, r$.
This follows from the cohomological exact sequence
\[H^{d_j}(X,F_j;\kk) \to H^{d_j}(X;\kk) \to H^{d_j}(F_j;\kk),\]
since $d_j>0$ and the second map vanishes as a consequence of the fact that the inclusion $F_j\to X$ is null-homotopic. 

Assuming that $N>0$ is large enough we may use the lifting property of Theorem \ref{lslift},
applied to compact $K=\bar U$, to find a lift of the cohomology class $u\in H^q(X;a^\xi)$ to a relative cohomology class
$\tilde u\in H^q(X, F\cap U; a^\xi)$. 
Similarly, by Theorem \ref{lslift} applied to the compact $\bar W$,
we may lift the class $v\in H^{q'}(X;b^\xi)$ to a relative cohomology class
$\tilde v\in H^{q'}(X, F\cap W; a^\xi)$. 

Therefore the product (\ref{lscup1}) is obtained from a product
\[\tilde u\cup \tilde v \cup \tilde w_1\cup \dots \cup \tilde w_r\, \in H^\ast(X,X;(ab)^\xi) =0\]
(lying in the trivial group) by restricting onto $X$. 
Hence any cup-product (\ref{lscup1}) must vanish. $\Box$

{\bf Remark.} Lemma 6.6 in \cite{F1} is incorrect as stated. 
As a consequence Proposition 6.5 and Corollary 6.7 of \cite{F1}
are correct only with some extra conditions. 
The lifting property of 
Theorem \ref{lslift} of the present paper
replaces Corollary 6.7 of \cite{F1}. 

\subsection{Cohomological estimate using Massey products} 

Let $X$ be a finite CW-complex and let $\xi\in H^1(X;\Z)$. 

\begin{definition} A cohomology class $v\in H^q(X;\kk)$ is called 
{\it $\xi$-survivor} if vanishes the cup-product $v\cup \xi=0$ and vanish all higher Massey products of the form
\[\langle v, \underbrace{\xi, \xi, \dots, \xi}_{r \, \, \mbox{times}}\rangle \, \in \, H^{q+1}(X;\kk)\]
for any $r>1$.
\end{definition}

We refer to sections 5 and 9 of paper \cite{F1}, where these Massey products were described in detail. 

\begin{theorem}\label{lsmassey} 
Let $X$ be a finite CW-complex and let $\xi\in H^1(X;\Z)$ be indivisible. 
Assume that there exist cohomology classes
$w_j \in H^{d_j}(X;\kk)$ of positive degree $d_j>0$, where $j=1, \dots, r$, having a nonzero cup-product 
$$0\neq w_1\cup w_2\cup \dots\cup w_r\in H^\ast(X;\kk).$$ 
If among the classes $w_1, \dots, w_r$ at least two are $\xi$-survivors, then 
\[\cat(X,\xi) \ge r-1.\]
\end{theorem}
{\bf Proof.} We may assume in the proof 
that the coefficient field $\kk$ is algebraically closed; otherwise we could replace $\kk$ by
its algebraic closure. 

Suppose that two first cohomology classes $w_1$, $w_2$ are $\xi$-survivors. 
It is shown in section 9.5 of \cite{F1}, that one may deform $w_1$ and $w_2$ to cohomology classes
$w'_1\in H^{d_1}(X;a^{\xi})$ and $w'_2\in H^{d_2}(X;a^{-\xi})$, where $a\in \kk$ is generic, 
such that the cup-product
\[0\not= w'_1\cup w_2'\cup w_3\cup \dots \cup w_r\in H^\ast(X;\kk)\]
is still nonzero. Now 
Theorem \ref{lsmassey} follows from Theorem \ref{lscup}. $\Box$

\begin{example}{\rm  Let $T^n=S^1\times \dots \times S^1$ be $n$-dimensional torus. 
Fix a point $x=(x_1, \dots, x_n)\in T^n$ and for $i=1, \dots, n$ consider $(n-1)$-dimensional subtorus
$T_i^{n-1}\subset T^n$ consisting of points with $i$-th coordinate equal to $x_i$. 
Let $X=T^n \# (S^1\times S^{n-1})$ be obtained from $T^n$ by adding a handle of index $1$. In other
words, we remove from $T^n$ interiors of two small disjoint disks $D_1$ and $D_2$ and connect their boundaries
by a tube $S^{n-1}\times [0,1]$. We will assume that the disks $D_1$ and $D_2$ do not meet the subtori
$T_i^{n-1}$ for $i=1, \dots, n$. Let $\xi\in H^1(X;\Z)$ be the cohomology class Poincar\'e dual to $\partial D_1$.

Each torus $T_i\subset X$ has trivial normal bundle and hence 
determines a cohomology class $w_i\in H^1(X;\Z)$. It is clear that 
the cup-product 
$$w_1\cup\dots\cup w_n\in H^n(X;\Z)$$ 
is nonzero since $T_1, \dots , T_n$ are mutually transversal and their intersection
consists of one point. As in \S 5 of \cite{F1} one may show that the classes $w_i$ are all $\xi$-survivors.
Theorem \ref{lsmassey} applies and gives $\cat(X,\xi)\geq n-1$. By Theorem \ref{lsmain1} any
smooth closed 1-form in class $\xi$ has at least $n-1$ geometrically distinct zeros.  

This estimate is sharp and cannot be improved. Indeed, since $\dim X =n$
Theorem 3.2 of \cite{F1} claims that there always exists a closed 1-form $\omega$ in class $\xi$ with $n-1$ zeros.}
\end{example}

\begin{example}{\rm Let $X$ be bouquet $X=Y\vee S^1$ as in Example \ref{lsexample14}. We will assume 
that the cohomology class
$\xi\in H^1(X;\Z)$ is such that
$\xi|_Y=0$ and $\xi|_{S^1}\neq 0$ is the generator. 
We want to find the estimate given by Theorem \ref{lsmassey}
in this example. We have $H^q(X;\kk) =H^q(Y;\kk)$ for $q>1$ and $H^1(X;\kk) = H^1(Y;\kk)\oplus H^1(S^1;\kk)$.
The last summand is generated by class $\xi$. It has trivial cup-products and Massey products with all other classes.  
If cohomology classes $w_j\in H^{d_j}(Y;\kk)$, where $j=1, \dots, r$, 
are such that $d_j>0$ and 
$$0\not= w_1\cup\dots \cup w_r \in H^\ast(Y;\kk)$$
then Theorem \ref{lsmassey} applies and gives $\cat(X,\xi)\ge r-1$. In Example \ref{lsexample14} we have shown
that $\cat(X, \xi) =\cat(Y)-1$. Hence we obtain $\cat(Y)\geq r$, which is weaker by 1 than the well-known inequality
claiming that $\cat(Y) $ is greater than the cup-length of $Y$.

This example shows that our cohomological estimates may be slightly improved.}
\end{example}

\appendix
\section{\large Appendix: Closed 1-forms on topological spaces}\label{lsapp1}

Differential 1-forms are defined only for smooth manifolds. 
{\it Closed} 1-forms may be defined for general topological spaces, as we show in this appendix. 

{\it A continuous closed 1-form $\omega$ on a topological space} 
$X$ is defined as a collection $\{f_U\}_{U\in \mathcal U}$
of continuous real valued functions $f_U: U\to \R$, where $\mathcal U=\{U\}$ is an open cover of $X$,
such that for any pair $U, V\in \mathcal U$ the difference
\[f_U|_{U\cap V} - f_V|_{U\cap V}: U\cap V \to \R\]
is a locally constant function.  Another such collection $\{g_V\}_{V\in \mathcal V}$ 
(where $\mathcal V$ is another open cover of $X$) defines {\it an equivalent} closed 1-form 
if for any point $x\in X$ there is an open neighborhood $W$ such that  for some
open sets $U\in \mathcal U$ and $V\in \mathcal V$ containing $W$ the difference
$f_U|_W \, -\,  g_V|_W$
is locally constant. 
The set of all continuous closed 1-forms on $X$ is a real vector space. 

As an example consider an open cover $\mathcal U=\{X\}$ consisting of the whole space $X$. Then any continuous 
function $f: X\to \R$ defines a closed 1-form on $X$, which is denoted $df$. 

For two continuous functions $f, g: X\to \R$ holds $df =dg$ 
if and only if the difference $f-g: X\to \R$ is locally constant.

One may integrate continuous closed 1-forms along continuous paths.
Let $\omega$ be a continuous closed 1-form  on $X$ given 
by a collection of continuous real valued functions $\{f_U\}_{U\in \mathcal U}$
with respect to an open cover $\mathcal U$ of $X$. Let $\gamma: [0,1]\to X$ be a continuous path. The line integral
$\int_\gamma \omega$ is defined as follows. Find a subdivision $t_0=0<t_1<\dots < t_N=1$ of the interval
$[0,1]$ such that for any $i$ the image $\gamma[t_i,t_{i+1}]$ is contained in a single open set $U_i\in \mathcal U$.
 Then we define
\begin{eqnarray}
\int\limits_\gamma \omega \, =\, 
\sum\limits_{i=0}^{N-1} \, \, [f_{U_i}(\gamma(t_{i+1})) - f_{U_i}(\gamma(t_i))].\label{lsintegral}
\end{eqnarray}
The standard argument shows that the integral (\ref{lsintegral}) does not depend on the choice of the subdivision and the 
open cover $\mathcal U$. 

\begin{lemma} Let $\omega$ be a continuous closed 1-form on topological spaces $X$.
For any pair of continuous
paths $\gamma, \gamma': [0,1]\to X$ with the 
common beginning $\gamma(0)=\gamma'(0)$ and the common end points $\gamma(1)=\gamma'(1)$, holds
$$\int_\gamma \omega \, =\,  \int_{\gamma'}\omega,$$
provided that $\gamma$ and $\gamma'$ are homotopic relative to the boundary.
\end{lemma}

{\bf Proof}. It is standard. $\Box$

Any closed 1-form defines {\it homomorphism of periods}
\begin{eqnarray}
\pi_1(X, x_0) \to \R, \quad [\gamma]\mapsto \int_\gamma \omega \in \R\label{lsperiod}
\end{eqnarray}
given by integration of 1-form $\omega$ along closed loops $\gamma: [0,1]\to X$ with $\gamma(0)=x_0=\gamma(1)$. 

\begin{lemma}
The homomorphism of periods (\ref{lsperiod}) is a group homomorphism.
\end{lemma}

\begin{lemma} Let $X$ be a path connected topological space. 
A continuous closed 1-form $\omega$ on $X$ equals $df$ for a
continuous function $f: X\to \R$ if and only if $\omega$ defines a trivial homomorphism of periods (\ref{lsperiod}).
\end{lemma}
{\bf Proof.} If $\omega =df$ then for any path $\gamma$ in $X$ holds $\int_\gamma \omega = f(q)-f(p)$, where
$q=\gamma(1), \, \, p=\gamma(0)$. Hence $\int_\gamma \omega =0$ if $\gamma$ is a closed loop. 

Conversely, assume that the homomorphism of periods
(\ref{lsperiod}) is trivial. One defines a continuous function $f:X\to \R$ by
$$f(x) \, =\, \int_{x_0}^x \omega.$$
Here the integration is taken over an arbitrary path connecting $x_0$ to $x$.
Assume that $\omega$ is given by a collection of continuous functions $f_U: U\to \R$ with respect to an
open cover $\{U\}$ of $X$. Then for any two points $x, y$ lying in the same path-connected component of $U$,
\[f(y) -f(x) =\int_x^y \omega = f_U(y) -f_U(x).\]
This shows that the function $f-f_U$ is locally constant on $U$. Hence
$df =\omega$. $\Box$

Any continuous 
closed 1-form $\omega$ on a topological space $X$ defines a (singular) {\it cohomology class} $[\omega]\in H^1(X;\R)$.
It is defined by the homomorphism of periods (\ref{lsperiod}) viewed as an element of 
$\Hom(H_1(X);\R)=H^1(X;\R)$. As follows from the above Lemma, {\it two continuous
closed 1-forms $\omega$ and $\omega'$ on 
$X$ have the same cohomology class $[\omega]=[\omega']$
if and only if their difference
$\omega - \omega'$ equals $df$, where $f:X\to \R$ is a continuous function. }

Recall that a topological space $X$ is {\it homologically locally connected} 
if for every point $x\in X$ and a neighborhood $U$ 
of $x$ there exists a neighborhood $V$ of $x$ in $U$ such that $\tilde H_q(V)\to \tilde H_q(U)$ is trivial for all $q$.

\begin{lemma} Let $X$ be a paracompact Hausdorff homologically locally connected topological space. 
Then any singular cohomology class $\xi\in H^1(X;\R)$ 
may be realized by a continuous closed 1-form on $X$.
\end{lemma}
\begin{proof} Consider the following exact sequence of sheaves over $X$
\begin{eqnarray}
0\to \R_X \to C_X\to B_X\to 0.\label{lsexact}
\end{eqnarray}
Here $\R_X$ denotes the sheaf of locally constant functions, $C_X$ denotes the sheaf of real valued 
continuous functions, and
$B_X$ denotes the sheaf of germs of continuous functions modulo locally constant.
More precisely, $B_X$ is the sheaf corresponding to the presheaf $U\mapsto C_X(U)/\R_X(U)$. 
Comparing with our definition of continuous closed 1-form, we find that {\it the space of global sections of $B_X$,
\[ H^0(X;B_X),\]
coincides with the space of continuous closed 1-forms on $X$}. 

From (\ref{lsexact}), using that $C_X$ is a fine sheaf, we obtain an exact sequence
\[0\to H^0(X;\R)\to H^0(X;C_X) \stackrel d\to H^0(X;B_X) \stackrel {[\, \, ]}\to H^1(X;\R_X)\to 0.\]
Here $H^0(X;C_X)=C(X)$ is the set of all continuous functions on $X$, and the map $d$ acts by 
assigning to a continuous function $f: X\to \R$ the closed 1-form $df\in H^0(X;B_X)$. 
The group $H^1(X;\R_X)$ is the ${\rm \check C}$ech cohomology $\check H^1(X;\R)$ and the map 
$[\, \, ]$ assigns to a closed 1-form $\omega$ its $ \rm\check C$ech cohomology class $[\omega]\in \check H^1(X;\R)$. 
The natural map $\check H^1(X;\R) \to H^1(X;\R)$ is an isomorphism assuming that $X$ is paracompact, Hausdorff, 
and homologically locally connected, cf. \cite{S}. 

This implies our statement.
$\Box$
\end{proof}

As we have shown in the proof, the space of continuous closed 1-forms $H^0(X;B_X)$ on a connected topological
space $X$ can be described 
by the following short exact sequence
\begin{eqnarray}
0\to C(X)/\R\, \stackrel d\to\,  H^0(X;B_X) \stackrel{[\, \, \, ]}{\, \to\,  }\check H^1(X;\R)\to 0.
\end{eqnarray}

\section{\large Appendix: Gradient-Convex Neighborhoods}\label{lsapp2}

 Let $M$ be a smooth manifold and let
$f:M\to \R$ be a $C^2$-smooth function with isolated critical points. 
Let $v$ be a gradient-like vector field for $f$.
This means that $v(f)>0$ on the complement of the set of critical points, and $v$ coincides with the gradient of $f$
on an open neighborhood of the critical points with respect to a Riemannian metric.

We will denote by $M\times \R\to M$, $(m,t)\mapsto m\cdot t$, the flow of the field $v$; we will assume that
it is defined for all $t\in \R$. 

\begin{lemma} Any open neighborhood $V\subset M$ of a critical point $p\in M$, $df_p=0$, 
contains a compact neighborhood $U$ of $p$, such that

\noindent
(1) for any point $m\in M$ the set $J_m=\{t\in \R;m\cdot t\in U\}\subset \R$ is either empty, or 
a closed interval $[a_m, b_m]$, possibly degenerated to a point, i.e. with $a_m=b_m$,

\noindent
(2) 
the function  $\{m\in M; J_m\neq \emptyset\}\to \R$, where $m\mapsto a_m\in \R$, is continuous.
\end{lemma}

The proof below was essentially suggested by P. Millman.

\begin{proof} 
We will assume that $f(p)=0$, that $\bar V$ is compact, and that point $p$ is the only critical point of function 
$f$ in $\bar V$.
Let $V_0$ be an open neighborhood of $p$ with compact closure $\bar V_0\subset V$
and such that in $V_0$ there exist local coordinates $x_1, \dots, x_n$ and $v|_{V_0}$ is the gradient of $f$ with respect to 
a Riemannian metric $g_{ij}$.
Fix a smooth function
$\psi: M\to [0,1]$, such that $\psi|_{V_0}\equiv 0$ and $\psi|_{(M-V)}\equiv 1$.

We want to show that there is a constant $\lambda>0$ such that the derivative of the function 
$$\Psi = \langle \grad f, \grad f \rangle: M\to \R_+$$ 
along the gradient flow of $f$ satisfies  in $V_0$
the inequalities:
\begin{eqnarray}
-\lambda \Psi \, \le\,  \frac{d \Psi}{dt} \, \le\,  \lambda \Psi.\label{lsestim}
\end{eqnarray}
In the local coordinates 
$$\Psi = \sum_{ij} g^{ij}\cdot\frac{\partial f}{\partial x_i}\cdot \frac{\partial f}{\partial x_j}$$
and 
\begin{eqnarray*}
 \frac{d \Psi}{dt} \, =\,  
\sum_{ij} \langle \grad \, g^{ij}, \grad f\rangle\cdot\frac{\partial f}{\partial x_i}\cdot \frac{\partial f}{\partial x_j} +\\
+ 2 \, \langle \grad f, \sum_{ij} g^{ij}\cdot \frac{\partial f}{\partial x_j}\cdot \grad \frac{\partial f}{\partial x_i}\rangle.
\end{eqnarray*}
Both the first and the second terms in this sum can be viewed as symmetric bilinear forms in the partial derivatives of $f$ 
with continuous coefficients,
and hence our statement (\ref{lsestim}) follows. 

We will define now two smooth functions $F_+, F_-: M\to \R$ by
\begin{eqnarray}
F_\pm \, =\,  \pm 2 f \, + \, \lambda^{-1} \cdot \Psi \, +\, \psi.
\end{eqnarray}
In $V_0-\{p\}$ using (\ref{lsestim}) we have 
$$
\begin{array}{l}
\displaystyle{
\frac{dF_+}{dt} \, = \, 2\Psi + \lambda^{-1} \cdot\frac{d\Psi}{dt} \, \geq \Psi \, >0,} \\ \\
\displaystyle{\frac{dF_-}{dt} \, = \,- 2\Psi + \lambda^{-1}\cdot \frac{d\Psi}{dt} \, \leq - \Psi\, <0.}
\end{array}
$$
Hence we conclude that $F_+$ increases and $F_-$ decreases along the gradient flow of $f$ in $V_0-\{p\}$. 

We set 
\[g =\max\{F_+, F_-\}: M \to \R_+.\]
It is a continuous function with $g(p)=0$. 
For any number $0<c$ small enough, the set $U_c=\{x\in M; g(x) < c\}$ is an open neighborhood
of $p$, contained in $V_0$.  

Let $V_1$ be an open 
neighborhood of $p$ with $\bar V_1\subset V_0$.  
Let $\epsilon >0$ be such that any trajectory $\gamma(t)\in M$ of the flow of $v$ with
$\gamma(t_1)\in \bar V_1$, and $\gamma(t_2)\in \overline{M-V_0}$, where $t_1<t_2$, satisfies 
\begin{eqnarray}
f(\gamma(t_2)) - f(\gamma(t_1)) >\epsilon.\label{lsgain}
\end{eqnarray}

We will show that $\bar U_c$ satisfies the conditions of the Lemma assuming that
 $0<c<\epsilon$ and $c$ is small enough, so that
$U_c\subset V_1$. 
Since $\partial \bar U_c=g^{-1}(c)$, a trajectory $\gamma(t)=m\cdot t$ enters the set $\bar U_c$ at $t=a$ if and only if 
$F_-(\gamma(a))=c$
and $F_+(\gamma(a)) \leq c$. 
Moreover, if $F_-(\gamma(a))=c$ and $F_+(\gamma(a))=c$, the trajectory leaves $\bar U_c$ immediately
(i.e. the trajectory $\gamma(t)$ is tangent to the boundary $\partial U_c$), 
and if $F_-(\gamma(a))< c$, the trajectory penetrates the interior of $\bar U_c$. 

Similarly, if $t=b$ is such that
$F_-(\gamma(b))\le c$ and $F_+(\gamma(b))=c$, the trajectory $\gamma(t)$ leaves the set $\bar U_c$ for $t>b$.
We know that, 
while $\gamma(t)$ stays in $V_1$, the function
 $F_+$ increases and so the trajectory remains away from $\bar U_c$.
May this trajectory come back to $\bar U_c$ for some large time $t=c>b$? If this happens, then $f(\gamma(c))>\epsilon$
(because of (\ref{lsgain}) and using $f(\gamma(b))\geq 0$) and hence $g(\gamma(c)) > \epsilon$. 

This proves that under our assumptions on $c$
the set $\{t;\gamma(t)\in \bar U_c\}$ coincides with the interval $[a,b]$. 

Now we are left to prove statement (2) of the Lemma. 
Let $A\subset M$ denote the set of points $m\in M$, 
such that $m\cdot t$ belongs to $V_1$ for some $t$. 
Since $\displaystyle{\frac{dF_-}{dt}}<0$, the equation
 $F_-(m\cdot a_m) =c$ defines a continuous function of $m\in A$. Similarly, the equation
 $F_+(m\cdot b_m) =c$ defines a continuous function $A\to \R$, where $m\mapsto b_m$.
The set $\{m\in M; J_m\not= \emptyset\}$ equals $\{m\in A; a_m\leq b_m\}$.  This implies our claim (2).
$\Box$
\end{proof} 

\bibliographystyle{amsalpha}

\vskip 2cm 

Address: 

Michael Farber,

Department of Mathematics, 

Tel Aviv University, Ramat Aviv 69978, Israel

farber@math.tau.ac.il

\end{document}